\documentclass[12pt]
{amsart}
\usepackage{latexsym,amscd,amssymb,comment}  
\makeatletter
\theoremstyle{plain}
  \newtheorem{thm}{Theorem}[section]
  \newtheorem{prop}[thm]{Proposition}
  \newtheorem{lem}[thm]{Lemma}
  \newtheorem{cor}[thm]{Corollary}
\theoremstyle{definition}
  \newtheorem{dfn}[thm]{Definition}
  \newtheorem{exmp}[thm]{Example}
  
\theoremstyle{remark}
  \newtheorem{rem}[thm]{Remark}
\numberwithin{equation}{section}
\def\ld{\operatorname{ld}}
\def\lin{\operatorname{lin}}
\def\SqR{\operatorname{Sq}(R)}
\def\SqS{\operatorname{Sq}(S)}
\def\SqE{\operatorname{Sq}(E)}
\def\<{{\langle}}
\def\>{{\rangle}}
\def\iff{\Longleftrightarrow}
\def\cmp{{\sf c}}
\def\Id{\operatorname{Id}}
\def\bA{{\bf A}}
\def\ba{{\bf a}}
\def\bb{{\bf b}}
\def\bc{{\bf c}}
\def\be{{\bf e}}

\def\bx{{\bf x}}
\def\00{{\bf 0}}
\def\b1{{\bf 1}}
\def\bP{ {\mathbf P}}
\def\cL{{\mathbf L}}
\def\R{{\mathbb R}}

\def\Z{{\mathbb Z}}
\def\NN{{\mathbb N}}
\def\rH{\tilde{H}}

\def\m{{\mathfrak m}}
\def\DE{{\bf D}_E}
\def\P{P_\bullet}
\def\I{I^\bullet}

\def\M{M^\bullet}

\def\cS{{\mathcal S}}
\def\cE{{\mathcal E}}
\def\Ext{\operatorname{Ext}}
\def\Hom{\operatorname{Hom}}
\def\depth{\operatorname{depth}}

\def\ann{\operatorname{ann}}
\def\pd{\operatorname{proj.dim}}
\def\chara{\operatorname{char}}
\def\RHom{{\rm R}{\Hom}}
\def\grR{\operatorname{*mod}R}

\def\grE{\operatorname{*mod}E}
\def\grS{\operatorname{*mod}S}
\def\can{\omega_R}
\def\wS{\omega_S}
\def\Dcom{\mathcal{D}^\bullet}

\def\const{\underline{K}}

\def\cExt{{\mathcal Ext}}

\def\Supp{\operatorname{supp}}
\def\nat{\operatorname{nat}} 
\def\Du{{\mathbf D}}
\def\relint{\operatorname{rel-int}}
\def\ldirr{\operatorname{ld.irr}}

\let\d=\Delta
\def\verset{\left[ n\right]}

\def\ubull#1{{#1}^{\bullet}}
\def\lbull#1{{#1}_{\bullet}}

\def\betti#1#2{\beta _{#1,#2}}

\def\PP{\mathbb{P}^2\mathbb{R}}

\def\tensor#1#2{\otimes ^{#1}_{#2}}
\def\dsum#1#2{\bigoplus ^{#1}_{#2}}
\def\ep#1{( #1)} 
\def\bra#1{\bigl[ #1\bigr]}  
\def\ebra#1{[ #1]} 
\def\sbra#1{\{ #1\}}  
\def\angle#1{\< #1\>}
\def\indeg{\operatorname{indeg}}

\def\ver{\operatorname{ver}}
\def\lk#1#2{\operatorname{lk}_{#1}{#2}}
\def\sk#1#2{{#1}^{\left({#2}\right)}}
\def\dual#1{{#1}^{\vee}}

\def\sign#1#2{\alpha \left( #1, #2\right)}
\def\redcomp#1#2{\tilde{\mathcal{C}}_{#1}\ep{#2;K}}
\def\redcocomp#1#2{\tilde{\mathcal{C}}^{#1}\ep{#2;K}}
\def\koszul{\mathcal{K}}

\def\tor#1#2#3#4{\operatorname{Tor}_{#1}^{#2}\ep{#3,#4}}

\def\redhom#1#2#3{\tilde{H}_{#1}\ep{#2;#3}}
\def\redcohom#1#2#3{\tilde{H}^{#1}\ep{#2;#3}}
\def\H#1#2{H_{#1}\ep{#2}}

\def\symfacering#1{K\ebra{#1}}
\def\exfacering#1{K\angle{#1}}
\def\stideali#1{I_{#1}}  
\def\stidealj#1{J_{#1}}  
\def\eref#1{(\ref{eq:#1})}
\title{Linearity Defects of Face Rings}
\dedicatory{Dedicated to Professor J\"urgen Herzog on his 65th birthday}
\author{Ryota Okazaki}
\address{Department of Mathematics, 
Graduate School of Science, Osaka University, Toyonaka, Osaka 
560-0043, Japan}
\email{smv679or@ecs.cmc.osaka-u.ac.jp}
\author{Kohji Yanagawa}
\address{Department of Mathematics, 
Graduate School of Science, Osaka University, Toyonaka, Osaka 
560-0043, Japan}
\email{yanagawa@math.sci.osaka-u.ac.jp}
\keywords{Stanley-Reisner ring, exterior face ring, linearity defect, 
weakly Koszul module, componentwise linear, sequentially Cohen-Macaulay,  
squarefree module}  
\begin{document}
\begin{abstract}
Let $S = K[x_1, \ldots, x_n ]$ be a polynomial ring over a field $K$, and 
$E = \bigwedge \< y_1, \ldots, y_n \>$ an exterior algebra. 
The {\it linearity defect} $\ld_E(N)$
of a finitely generated graded $E$-module $N$ 
measures how far $N$ departs from ``componentwise linear". 
It is known that $\ld_E(N) < \infty$ for all $N$. But  
the value can be arbitrary large, while 
the similar invariant $\ld_S(M)$ for an $S$-module $M$ 
is always at most $n$. We will show that if $I_\Delta$ (resp. $J_\Delta$) 
is the squarefree monomial ideal of $S$ (resp. $E$) corresponding 
to a simplicial complex $\Delta \subset 2^{\{1, \ldots, n \}}$, 
then $\ld_E(E/J_\Delta) = \ld_S(S/I_\Delta)$. Moreover, 
except some extremal cases, 
$\ld_E(E/J_\Delta)$ is a topological invariant of the geometric realization 
$|\Delta^\vee|$ of the Alexander dual $\Delta^\vee$ of $\Delta$. 
We also show that, 
when $n \geq 4$, $\ld_E(E/J_\Delta) = n-2$ (this is the largest possible 
value) if and only if $\Delta$ is an $n$-gon. 
\end{abstract}

\maketitle

\section{Introduction}
Let $A = \bigoplus_{i \in \NN} A_i$ be a graded 
(not necessarily commutative) noetherian algebra over a field 
$K \, (\cong A_0)$. 
Let $M$ be a finitely generated graded left $A$-module, 
and $\P$ its minimal free resolution.  
Eisenbud et al. \cite{EFS} defined 
the {\it linear part} $\lin(\P)$ of $\P$,  
which is the complex obtained by erasing all terms of degree 
$\geq 2$ from the matrices representing the differential maps 
of $\P$ (hence $\lin(\P)_i = P_i$ for all $i$).   
Following Herzog and Iyengar \cite{HI},  
we call  $\ld_A(M) = \sup\{ \, i \mid H_i(\lin(\P)) \ne 0 \, \}$
the {\it linearity defect} of $M$. 
This invariant and related concepts have been studied by 
several authors (e.g., \cite{EFS, HI, MZ, R02, Y7}). 
We say a finitely generated graded $A$-module $M$ is 
{\it componentwise linear} (or, 
{\it (weakly) Koszul} in some literature) if $M_{\< i \>}$ has a linear 
free resolution for all $i$. Here $M_{\< i \>}$ is  
the submodule of $M$ generated by its degree $i$ part $M_i$. 
Then we have $$\ld_A(M) = \min \{ \, i \mid 
\text{the $i^{\rm th}$ syzygy of $M$ is componentwise linear} \, \}.$$

For this invariant, a remarkable result holds 
over an exterior algebra $E = \bigwedge \< y_1, \ldots, y_n \>$. 
In \cite[Theorem~3.1]{EFS}, 
Eisenbud et al. showed that any finitely generated graded $E$-module 
$N$ satisfies $\ld_E(N) < \infty$ while $\pd_E(N) = \infty$ in most cases. 
(We also remark that Martinez-Villa and 
Zacharia~\cite{MZ} proved the same result for many selfinjective 
Koszul algebras). 
If $n \geq 2$, then we have 
$\sup \{ \, \ld_E(N) \mid  \text{$N$ a finitely generated graded 
$E$-module} \, \} = \infty$. 
But Herzog and R\"omer proved that if $J \subset E$ is a 
{\it monomial} ideal then $\ld_E(E/J) \leq n-1$ (c.f. \cite{R02}).   

A monomial ideal of $E = \bigwedge \< y_1, \ldots, y_n \>$ is always 
of the form $J_\Delta := ( \, \prod_{i \in F} y_i \mid  
F \not \in \Delta \, )$ for a simplicial complex 
$\Delta \subset 2^{\{1, \ldots, n \}}$. 
Similarly, we have the {\it Stanley-Reisner ideal}  
$I_\Delta := ( \, \prod_{i \in F} x_i \mid  F \not \in \Delta \, )$ 
of a polynomial ring $S=K[x_1, \ldots, x_n ]$. 
In this paper, we will show the following. 

\begin{thm}
With the above notation, we have $\ld_E(E/J_\Delta) = 
\ld_S (S/I_\Delta)$. Moreover, 
if $\ld_E(E/J_\Delta) > 0$ (equivalently, 
$\Delta \ne 2^{T}$ for any $T \subset [n]$),
then $\ld_E(E/J_\Delta)$ is a 
topological invariant of the geometric realization 
$|\Delta^\vee|$ of the Alexander dual $\Delta^\vee$.
(But $\ld(E/J_\Delta)$ may depend on $\chara(K)$.)
\end{thm}

By virtue of the above theorem, we can put 
$\ld(\Delta):= \ld_E(E/J_\Delta) = \ld_S (S/I_\Delta)$. 
If we set $d:= \min \{ \,i \mid [I_\Delta]_i \ne 0 \, \} = 
\min \{ \,i \mid [J_\Delta]_i \ne 0 \, \}$, then 
$\ld(\Delta) \leq \max \{1, n- d\}$. 
But, if $d=1$ (i.e., 
$\{ i \} \not \in \Delta$ for some $1 \leq i \leq n$), 
then $\ld(\Delta) \leq \max \{1, n-3\}$. Hence, 
if $n \geq 3$, we have $\ld(\Delta) \leq n-2$ for all $\Delta$.

\begin{thm}
Assume that $n \geq 4$. Then $\ld(\Delta) = n-2$ 
if and only if $\Delta$ is an $n$-gon. 
\end{thm}

While we treat $S$ and $E$ in most part of the paper, some results 
on $S$ can be generalized to a normal semigroup ring, 
and this generalization makes the topological meaning 
of $\ld(\Delta)$ clear. 
So \S2 concerns a normal semigroup ring. But, in this case, we use 
an irreducible resolution (something analogous 
to an injective resolution), not a projective resolution.

\section{Linearity Defects for Irreducible Resolutions}
Let $C \subset \Z^n \subset \R^n$ be an affine semigroup 
(i.e., $C$ is a finitely generated additive submonoid of $\Z^n$), 
and $R := K[\bx^\bc \mid \bc \in C] \subset K[x_1^{\pm 1}, \ldots, 
x_n^{\pm 1}]$ the semigroup ring of $C$ over the field $K$. 
Here $\bx^\bc$ for $\bc = (c_1, \ldots, c_n) \in C$  
denotes the monomial $\prod_{i=1}^n x_i^{c_i}$.  
Let $\bP := \R_{\geq 0} C \subset \R^n$ be the polyhedral cone 
spanned by $C$.  
We always assume that $\Z C = \Z^n$, $\Z^n \cap \bP = C$ 
and $C \cap (-C) = \{ 0 \}$. 
Thus $R$ is a normal Cohen-Macaulay integral domain of dimension $n$ 
with a maximal ideal
$\m := (\bx^\bc \mid 0 \ne \bc \in C)$. 

Clearly, $R = \bigoplus_{\bc \in C} K \bx^\bc$ is a $\Z^n$-graded ring. 
We say a $\Z^n$-graded ideal of $R$ is a {\it monomial ideal}. 
Let $\grR$ be the category of finitely generated $\Z^n$-graded 
$R$-modules and degree preserving $R$-homomorphisms. 
As usual, for $M \in \grR$ and $\ba \in \Z^n$, $M_\ba$ 
denotes the degree $\ba$ component of $M$, and $M(\ba)$ denotes the 
shifted module of $M$ with $M(\ba)_\bb = M_{\ba + \bb}$. 

Let $\cL$ be the set of non-empty faces of the polyhedral cone $\bP$.  
Note that $\{ 0\}$ and $\bP$ itself belong to $\cL$. For $F \in \cL$, 
$P_F := ( \, \bx^\bc \mid \bc \in C \setminus  F \, )$ is a prime ideal 
of $R$.  Conversely, any monomial prime ideal is of the 
form $P_F$ for some $F \in \cL$. Note that $P_{\{ 0\}} = \m$ 
and $P_\bP = (0)$. Set $K[F] := R/P_F \cong K[ \, \bx^\bc \mid 
\bc \in C \cap F]$ for $F \in \cL$. 
The Krull dimension of $K[F]$ equals the dimension $\dim F$ of 
the polyhedral cone $F$. 

For a point $u \in \bP$,  we always have a unique face 
$F \in \cL$ whose relative interior contains $u$. 
Here we denote $s(u) = F$.  

\begin{dfn}[\cite{Y2}]\label{sq}
We say a module $M \in \grR$ 
is {\it squarefree}, if it is $C$-graded 
(i.e., $M_\ba = 0$ for all $\ba \not \in C$),   
and the multiplication map 
$M_\ba \ni y \mapsto \bx^\bb y 
\in M_{\ba + \bb}$ is bijective for all 
$\ba, \bb \in C$ with $s(\ba+\bb) = s(\ba)$. 
\end{dfn} 

For a monomial ideal $I$, $R/I$ is a squarefree $R$-module 
if and only if $I$ is a radical ideal (i.e., $\sqrt{I} = I$). 
Regarding $\cL$ as a partially ordered set by inclusion, 
we say $\Delta \subset \cL$ is an {\it order ideal}, if 
$\Delta \ni F \supset F' \in \cL$ implies $F' \in \Delta$. 
If $\Delta$ is an order ideal, then 
$I_\Delta := ( \, \bx^\bc \mid \bc \in C, \,  
s(\bc) \not \in \Delta  \, ) \subset R$ 
is a radical ideal. 
Conversely, any radical monomial ideal is of the 
form $I_\Delta$ for some $\Delta$. 
Set $K[\Delta] := R/I_\Delta$. Clearly, 
$$
K[\Delta]_\ba \cong 
\begin{cases}
K & \text{if $\ba \in C$ and $s(\ba) \in \Delta$,}\\
0 & \text{otherwise.}
\end{cases}
$$
In particular, if $\Delta = \cL$ (resp. $\Delta = \{ \, \{ 0 \} \, \}$), 
then $I_\Delta = 0$ (resp. $I_\Delta = \m$) 
and $K[\Delta] = R$ (resp. $K[\Delta] = K$).  
When $R$ is a polynomial ring, $K[\Delta]$ is nothing else than 
the Stanley-Reisner ring of a simplicial complex $\Delta$. 
(If $R$ is a polynomial ring, then the partially ordered set $\cL$ 
is isomorphic to the power set $2^{\{1, \ldots, n\}}$, 
and $\Delta$ can be seen as a simplicial complex.) 

For each $F \in \cL$, take some $\bc(F) \in C \cap \relint(F)$ 
(i.e., $s(\bc(F)) =F$). For a squarefree $R$-module $M$ and 
$F, G \in \cL$ with $G \supset F$, \cite[Theorem~3.3]{Y2} gives a 
$K$-linear map $\varphi^M_{G, F}: M_{\bc(F)} \to 
M_{\bc(G)}$.  
They satisfy $\varphi^M_{F,F} = \Id$ and 
$\varphi^M_{H, G} \circ \varphi^M_{G, F} = \varphi^M_{H,F}$ for all 
$H \supset G \supset F$. We have $M_\bc \cong M_{\bc'}$ for 
$\bc, \bc' \in C$ with $s(\bc) = s(\bc')$. Under these isomorphisms, 
the maps  $\varphi^M_{G, F}$ do not depend on the particular 
choice of $\bc(F)$'s. 

Let $\SqR$ be the full subcategory of $\grR$ 
consisting of squarefree modules. As shown in \cite{Y2}, 
$\SqR$ is an abelian category with enough injectives. 
For an indecomposable squarefree module $M$, 
it is injective in $\SqR$ if and only if  
$M \cong K[F]$ for some $F \in \cL$.  
Each $M \in \SqR$ has a minimal injective resolution in $\SqR$, 
and we call it  a {\it minimal irreducible resolution} 
(see \cite{Y8} for further information). 
A minimal irreducible resolution is unique up to isomorphism, 
and its length is at most $n$. 

Let $\can$ be the $\Z^n$-graded canonical module of $R$. 
It is well-known that $\can$ is isomorphic to the radical 
monomial ideal $(\, \bx^\bc \mid \bc \in C, \, s(c) = \bP \, )$. 
Since we have $\Ext^i_R(\M,\can) \in \SqR$ for all $\M \in \SqR$, 
$\Du(-) := \RHom_R(-, \can)$ gives a duality functor from 
the derived category $D^b(\SqR) \, (\cong D^b_{\SqR} (\grR))$ to itself. 

In the sequel, for a $K$-vector space $V$, 
$V^*$ denotes its dual space.  But, even if $V = M_\ba$  
for some $M \in \grR$ and $\ba \in \Z^n$, 
we set the degree of $V^*$ to be 0.  

\begin{lem}[{\cite[Lemma~3.8]{Y8}}]\label{explicit form}
If $M \in \SqR$, then $\Du(M)$ is quasi-isomorphic to the complex 
$D^{\bullet} : 0 \to D^0 \to D^1 \to \cdots \to D^n \to 0$ with 
$$D^i = \bigoplus_{\substack{F \in \cL \\ \dim F = n-i}} 
(M_{\bc(F)})^* \otimes_K K[F].$$
Here the differential is the sum of the maps 
$$(\pm\varphi_{F,F'}^M)^* \otimes \nat
: (M_{\bc(F)})^* \otimes_K K[F] \to (M_{\bc(F')})^* \otimes_K 
K[F']$$
for  $F,F' \in \cL$ with $F \supset F'$ and $\dim F = \dim F'+1$, 
and $\nat$ denotes the natural surjection $K[F] \to K[F']$. 
We can also describe $\Du(\M)$ for a complex $\M \in D^b(\SqR)$ 
in a similar way. 
\end{lem}

\noindent{\bf Convention.} In the sequel, as an explicit complex, 
$\Du(\M)$ for $\M \in D^b(\SqR)$ means the complex described in 
Lemma~\ref{explicit form}. 

\medskip

Since $\Du \circ \Du \cong \Id_{D^b(\SqR)}$, $\Du \circ \Du (M)$ 
is an irreducible resolution of $M$, but it is far from being minimal. 
Let $(\I, \partial^\bullet)$ be a minimal irreducible resolution of $M$. 
For each $i \in \NN$ and $F \in \cL$, we have a natural number 
$\nu_i(F, M)$ such that 
$$I^i \cong \bigoplus_{F \in \cL} K[F]^{\nu_i(F,M)}.$$
Since $\I$ is minimal, $z \in K[F] \subset I^i$ with $\dim F = d$ 
is sent to 
$$\partial^i(z) \in \bigoplus_{\substack{G \in \cL \\ 
\dim G < d}} K[G]^{\nu_{i+1}(G,M)} \subset I^{i+1}.$$
The above observation on $\Du \circ \Du(M)$ gives the formula 
(\cite[Theorem~4.15]{Y2})
$$\nu_i(F,M) = \dim_K [\Ext^{n-i-\dim F}_R (M, \can)]_{\bc(F)}.$$ 

For each $l \in \NN$ with $0 \leq l \leq n$, 
we define the $l$-{\it linear strand} 
$\lin_l(\I)$ of $\I$ as follows: 
The term $\lin_l(\I)^i$ of cohomological degree $i$  
is $$\bigoplus_{\dim F = l-i} K[F]^{\nu_i(F,M)},$$
which is a direct summand of $I^i$, 
and the differential $\lin_l(\I)^i \to \lin_l(\I)^{i+1}$ 
is the corresponding component of the differential 
$\partial^i: I^i \to I^{i+1}$ of $\I$.  
By the minimality of $\I$, we can see that
$\lin _l \ep{\I}$ are cochain complexes.
Set $\lin(\I) := \bigoplus_{0 \leq l \leq n} \lin_l(\I)$. 
Then we have the following. 
For a complex $\M$ and an integer $p$, let $\M[p]$ be the $p^{\rm th}$ 
translation of $\M$. That is, $\M[p]$ is a complex with $M^i[p] = M^{i+p}$.

\begin{thm}[{\cite[Theorem~3.9]{Y8}}]\label{lin(M)}
With the above notation, we have  
$$\lin_l(\I) \cong \Du(\Ext_R^{n-l}(M,\can))[n-l].$$
Hence 
$$\lin(\I) \cong \bigoplus_{i \in \Z} \Du(\Ext_R^i(M,\can))[i].$$
\end{thm}

\begin{dfn}
Let $\I$ be a minimal irreducible resolution of $M \in \SqR$. 
We call
$\max \{ \, i \mid H^i(\lin(\I)) \ne 0 \, \}$ 
the {\it linearity defect of the minimal irreducible resolution} of $M$, 
and denote it by $\ldirr_R(M)$. 
\end{dfn}

\begin{cor}\label{ld for irr}
With the above notation, we have  
$$\max\{ \, i \mid H^i(\lin_l(\I)) \ne 0 \, \} 
= l  - \depth_R ( \, \Ext_R^{n-l}( M, \can ) \,),$$  
and hence 
$$\ldirr_R(M) = \max \{ \, i - \depth_R 
( \, \Ext_R^{n-i}( M, \can ) \, )  \mid  0 \leq i \leq n  \, \}.$$ 
Here we set the depth of the 0 module to be $+ \infty$.
\end{cor}

\begin{proof}
By Theorem~\ref{lin(M)}, we have $H^i(\lin_l(\I)) = 
\Ext^{i+l}_R(\Ext^l_R(M,\can), \can)$. 
Since $\depth_R N = \min \{ \, i \mid \Ext^{n-i}_R(N, \can) \ne 0 \, \}$ 
for a finitely generated graded $R$-module $N$, the assertion follows. 
\end{proof}

\begin{dfn}[Stanley \cite{St}] Let $M \in \grR$. 
We say $M$ is {\it sequentially Cohen-Macaulay} if there is a 
finite filtration 
$$0 = M_0 \subset M_1 \subset \cdots \subset M_r = M$$ 
of $M$ by graded submodules $M_i$ satisfying the following conditions.  
\begin{itemize}
\item[(a)] Each quotient $M_i/M_{i-1}$ is Cohen-Macaulay.
\item[(b)] $\dim(M_i/M_{i-1}) < \dim (M_{i+1}/M_i)$ for all $i$. 
\end{itemize}
\end{dfn} 

Remark that the notion of sequentially Cohen-Macaulay module is 
also studied under the name of a ``Cohen-Macaulay filtered module" (\cite{Sc}). 

Sequentially Cohen-Macaulay property is getting important in the 
theory of Stanley-Reisner rings. 
It is known that $M \in \grR$ is sequentially Cohen-Macaulay 
if and only if $\Ext^{n-i}_R(M,\can)$ is a zero module or a Cohen-Macaulay module 
of dimension $i$ for all $i$ (c.f. \cite[III. Theorem~2.11]{St}). 
Let us go back to Corollary~\ref{ld for irr}. 
If $N :=  \Ext_R^{n-i}( M, \can) \ne 0$,  
then $\depth_R N \leq \dim_R N \leq i$. 
Hence $\depth_R N = i$ if and only if 
$N$ is a Cohen-Macaulay module of dimension $i$. 
Thus, as stated in \cite[Corollary~3.11]{Y8}, 
$\ldirr_R(M)=0$ if and only if $M$ is sequentially Cohen-Macaulay.   

Let $\I: 0 \to I^0 \stackrel{\partial^0}{\to} I^1 
\stackrel{\partial^1}{\to} I^2 \to \cdots$ 
be an irreducible resolution of $M \in \SqR$. 
Then it is easy to see that 
$\ker(\partial^i)$ is sequentially Cohen-Macaulay if and only if 
$i \geq \ldirr_R(M)$. In particular, 
$$\ldirr_R(M) = \min\{ \, i \mid \text{$\ker(\partial^i)$ is 
sequentially Cohen-Macaulay}\}.$$

\bigskip

We have a hyperplane $H \subset \R^n$ such that $B:= H \cap \bP$ 
is an $(n-1)$-dimensional polytope. 
Clearly, $B$ is homeomorphic to a closed ball of dimension $n-1$. 
For a face $F \in \cL$, set $|F|$ to be the relative interior of $F \cap H$. 
If $\Delta \subset \cL$ is an order ideal,  
then $|\Delta|  := \bigcup_{F \in \Delta} |F|$ is a closed subset of $B$, and 
$\bigcup_{F \in \Delta} |F|$ is a {\it regular cell decomposition}
(c.f. \cite[\S 6.2]{BH}) of $|\Delta|$.  
Up to homeomorphism, (the regular cell decomposition of) $|\Delta|$ does 
not depend on the particular choice of the hyperplane $H$. 
The dimension  $\dim |\Delta|$ of $|\Delta|$ is given by 
$ \max \{ \, \dim |F| \mid F \in \Delta \, \}$. Here $\dim |F|$ denotes 
the dimension of $|F|$ as a cell (we set $\dim \emptyset = -1$), 
that is, $\dim |F| = \dim F -1 = \dim K[F]-1$. 
Hence we have $\dim K[\Delta] = \dim |\Delta| +1$. 

If $F \in \Delta$, then $U_F := \bigcup_{F' \supset F} |F'|$ 
is an open set of $B$. 
Note that $\{ \, U_F \mid \{ 0 \} \ne F \in \cL \, \}$ 
is an open covering of $B$. 
In \cite{Y6}, from $M \in \SqR$, we constructed a sheaf $M^+$ on $B$. 
More precisely, the assignment 
$$\Gamma(U_F, M^+) = M_{\bc(F)}$$ for each $F \ne \{ 0  \}$ and 
the map 
$$\varphi_{F,F'}^M:\Gamma(U_{F'},M^+) = M_{\bc(F')} \to M_{\bc(F)} = 
\Gamma(U_F,M^+)$$ for $F, F' \ne \{  0 \}$ with 
$F \supset F'$ (equivalently, $U_{F'} \supset U_{F}$) defines 
a sheaf. Note that $M^+$ is a {\it constructible sheaf} 
with respect to the cell decomposition $B=\bigcup_{F \in \cL} |F|$. 
In fact, for all $\{ 0 \} \ne F \in \cL$, 
the restriction $M^+|_{|F|}$ of $M^+$ to $|F| \subset B$ 
is a constant sheaf with coefficients in $M_{\bc(F)}$.  
Note that $M_\00$ is ``irrelevant" to $M^+$, where 
$\00$ denotes $(0,0, \ldots, 0) \in \Z^n$. 

It is easy to see that $K[\Delta]^+ \cong j_* \const_{|\Delta|}$, 
where $\const_{|\Delta|}$ is the constant sheaf on $|\Delta|$ 
with coefficients in $K$, 
and $j$ denotes the embedding map $|\Delta| \hookrightarrow B$. 
Similarly, we have that $(\can)^+ \cong h_! \const_{B^\circ}$, 
where $\const_{B^\circ}$ is the constant sheaf on 
the relative interior $B^\circ$ of $B$, 
and $h$ denotes the embedding map $B^\circ \hookrightarrow B$. 
Note that $(\can)^+$ is the orientation sheaf of $B$ over $K$. 

\begin{thm}[{\cite[Theorem~3.3]{Y6}}]\label{Hoch}
For $M \in \SqR$, we have an isomorphism 
$$H^i(B; M^+) \cong [H_\m^{i+1}(M)]_{\00} \quad  \text{for all $i \geq 1$},$$ 
and an exact sequence 
$$0 \to [H_\m^0(M)]_\00 \to M_\00 \to H^0( B; M^+) \to [H_\m^1(M)]_\00 \to 0.$$ 
In particular, we have 
$[H_\m^{i+1}(K[\Delta])]_\00 \cong \rH^i(|\Delta| ; K)$ for all $i \geq 0$, 
where $\rH^i(|\Delta| ; K)$ denotes the $i^{ th}$ reduced cohomology 
of $|\Delta|$ with coefficients in $K$. 
\end{thm}

Let $\Delta \subset \cL$ be an order ideal and $X := |\Delta|$. 
Then $X$ admits Verdier's dualizing complex $\Dcom_X$, 
which is a complex of sheaves of $K$-vector spaces. 
For example, $\Dcom_B$ is quasi-isomorphic to $(\can)^+[n-1]$. 

\begin{thm}[{\cite[Theorem~4.2]{Y6}}]\label{Verdier}
With the above notation, 
if $\ann(M) \supset I_\Delta$ 
(equivalently, $\Supp (M^+) := \{ x \in B \mid (M^+)_x \ne 0 \} 
\subset X$), then we have 
$$\Supp (\Ext_R^i(M, \can)^+) \subset X \quad  \text{and} \quad  
\Ext_R^i(M, \can)^+ |_X \cong \cExt^{i-n+1}(M^+|_X, \Dcom_X).$$
\end{thm}
 
\begin{thm}\label{ldirr sheaf}
Let $M$ be a squarefree $R$-module with $M \ne 0$ and 
$[H_\m^1(M)]_\00 = 0$, and $X$ the closure of $\Supp (M^+)$. 
Then $\ldirr_R(M)$ only depends on 
the sheaf $M^+|_X$ (also independent from $R$). 
\end{thm}

\begin{proof}
We use Corollary~\ref{ld for irr}. 
In the notation there, the case when $i = 0$ is always 
unnecessary to check. 
Moreover, by the present assumption, we have 
$\depth_R (\, \Ext_R^{n-1}(M, \can) \, ) \geq 1$ (in fact, 
$\Ext_R^{n-1}(M, \can)$ 
is either the 0 module, or a 1-dimensional Cohen-Macaulay 
module). So we may assume that $i > 1$.    

Recall that 
$$\depth_R ( \, \Ext_R ^{n-i}(M, \can) \, ) = \min  \{ \, j 
\mid \Ext^{n-j}_R( \, \Ext_R^{n-i}(M, \can), \can \, ) \ne 0 \, \}.$$
By Theorem~\ref{Verdier}, 
$[\Ext^{n-j}_R( \, \Ext_R^{n-i}(M, \can), \can \, )]_\ba$ can be determined by 
$M^+|_X$ for all $i,j$ and all $\ba \ne 0$.  If $j > 1$, then  
$[\Ext^{n-j}_R(\, \Ext_R^{n-i}(M, \can), \can \, )]_\00$ is isomorphic to 
\begin{eqnarray*}
[H_\m^j(\Ext_R^{n-i}(M, \can))]_\00^* &\cong& 
H^{j-1}( \, B; \,  \Ext_R^{n-i}(M, \can)^+ \, )^* \\
&\cong& 
H^{j-1}( \, X ; \,  \cExt^{-i-1}(M^+ |_X; \Dcom_X) \,)^*
\end{eqnarray*} 
(the first and the second isomorphisms follow from 
Theorem~\ref{Hoch} and Theorem~\ref{Verdier}, respectively),  
and determined by $M^+|_X$. 
So only $[\Ext^{n-j}_R(\Ext_R^{n-i}(M, \can), \can)]_\00$ for 
$j= 0,1$ remain. As above, they are isomorphic to 
$[H_\m^j(\Ext_R^{n-i}(M, \can))]_\00^*$. 
But, by \cite[Lemma~5.11]{Y8}, 
we can compute $[H_\m^j(\Ext_R^{n-i}(M, \can))]_\00$ 
for $i > 1$ and $j= 0,1$ from the sheaf  $M^+|_X$. 
So we are done. 
\end{proof}

\begin{thm}\label{topological}
For an order ideal $\Delta \subset \cL$ with $\Delta \ne \emptyset$, 
$\ldirr_R(K[\Delta])$ depends 
only on the topological space $|\Delta|$.  
\end{thm}

Note that $\ldirr_R(K[\Delta])$ may depend on $\chara(K)$. 
For example, if $|\Delta|$ is homeomorphic to a real projective plane, 
then $\ldirr_R(K[\Delta])=0$ if $\chara(K) \ne 2$, but 
$\ldirr_R(K[\Delta])=2$ if $\chara(K) =2$.  

Similarly, some other invariants and 
conditions (e.g., the Cohen-Macaulay property of $K[\Delta]$) 
studied in this paper depend on $\chara(K)$. 
But, since we fix the base field $K$, we always omit the phrase ``over $K$".

\begin{proof}
If $|\Delta|$ is not connected, then $[H_\m^1(K[\Delta])]_\00 \ne 0$ 
by Theorem~\ref{Hoch}, 
and we cannot use Theorem~\ref{ldirr sheaf} directly. 
But even in this case, $\depth_R (\,  \Ext^{n-i}_R(K[\Delta], \can) \, )$ 
can be computed for all $i \ne 1$ 
by the same way as in Theorem~\ref{ldirr sheaf}.
In particular, they only depend on $|\Delta|$. 
So the assertion follows from the next lemma. 
\end{proof}

\begin{lem}
We have $\depth_R  (\, \Ext^{n-1}_R(K[\Delta], \can) \, ) \in \{ 0,1, + \infty \}$, 
and 
$$\depth_R ( \, \Ext^{n-1}_R(K[\Delta], \can) \, )= 0 
\ \,   \text{if and only if} \ \,  
\text{$|\Delta'|$ is not connected.}$$ Here 
$\Delta' := \Delta \setminus \{ \, F \mid 
\text{$F$ is a maximal element of $\Delta$ and $\dim |F| = 0 $} \, \}.$
\end{lem}

\begin{proof}
Since $\dim_R \Ext^{n-1}_R(K[\Delta], \can) \leq 1$, 
the first statement is clear. If $\dim |\Delta| \leq 0$, then 
$|\Delta'| = \emptyset$ and $\depth_R (\, \Ext^{n-1}_R(K[\Delta], \can) \,) 
\geq 1$. So, to see the second statement, 
we may assume that $\dim |\Delta| > 1$. 
Set $J := I_{\Delta'}/I_\Delta$ to 
be an ideal of $K[\Delta]$. Note that either $J$ is a 
1-dimensional Cohen-Macaulay module or $J=0$. 
From the short exact sequence 
$0 \to J \to K[\Delta] \to K[\Delta'] \to 0$, we have an exact 
sequence $$0 \to \Ext_R^{n-1}(K[\Delta'],\can) \to 
\Ext^{n-1}_R(K[\Delta],\can) \to \Ext_R^{n-1}(J,\can)\to 0.$$
Since $\Ext_R^{n-1}(J,\can)$ has positive depth, 
$\depth_R (\, \Ext_R^{n-1}(K[\Delta'],\can) \, )=0$ if and only if 
$\depth_R (\, \Ext^{n-1}_R(K[\Delta],\can) \, ) =0$. 
But, since $K[\Delta']$ does not have 1-dimensional associated primes, 
$\Ext^{n-1}_R(K[\Delta'],\can)$ is an artinian module. 
Hence we have the following. 
\begin{eqnarray*}
\depth_R (\, \Ext^{n-1}_R(K[\Delta'],\can) \, ) =0  &\iff&  
[\Ext^{n-1}_R(K[\Delta'],\can)]_\00 \ne 0 \\ &\iff& 
[H_\m^1(K[\Delta'])]_\00 = \rH^0(|\Delta'|;K) \ne 0 \\ &\iff&  
\text{$|\Delta'|$ is not connected.}
\end{eqnarray*}
\end{proof}

\section{Linearity Defects of Symmetric and Exterior Face Rings}
Let $S := K[x_1, \ldots, x_n]$ be a polynomial ring,  
and consider its natural $\Z^n$-grading. 
Since $S = K[\NN^n]$ is a normal semigroup ring, we can use the notation and 
the results in the previous section. 

Now we introduce some conventions which 
are compatible with the previous notation. 
Let $\be_i := (0, \ldots, 0,1,0, \ldots 0) \in \R^n$ be the 
$i^{\rm th}$ unit vector, and $\bP$ the cone 
spanned by $\be_1, \ldots, \be_n$.  We identify a face $F$ of 
$\bP$ with the subset $\{ \, i \mid \be_i \in F \, \}$ of 
$[n] := \{ 1,2, \ldots, n \}$. 
Hence the set $\cL$ of nonempty faces of $\bP$ can be identified 
with the power set $2^{[n]}$ of $[n]$. 
We say $\ba = (a_1, \ldots, a_n)\in \NN^n$ is 
{\it squarefree}, if $a_i = 0,1$  for all $i$. 
A squarefree vector $\ba \in \NN^n$ will be identified 
with the subset $\{ \, i  \mid a_i = 1 \, \}$ of $[n]$. 
Recall that we took a vector $\bc(F) \in C$ for each $F \in \cL$ 
in the previous section. Here we assume that 
$\bc(F)$ is the squarefree vector corresponding to 
$F \in \cL \cong 2^{[n]}$. 
So, for a $\Z^n$-graded $S$-module $M$, we simply denote 
$M_{\bc(F)}$ by $M_F$. In the first principle, 
we regard $F$ as a subset of $[n]$, or a squarefree vector 
in $\NN^n$, rather than the corresponding face of $\bP$. 
For example, we write $P_F = (x_i \mid i \not \in F)$, 
$K[F] \cong K[x_i \mid i \in F]$. And $S(-F)$ denotes the rank 1 
free $S$-module $S(-\ba)$, where $\ba \in \NN^n$ is the squarefree 
vector corresponding to $F$. 

Squarefree $S$-modules are defined by the same way as Definition~\ref{sq}. 
Note that the free module $S(-\ba)$, $\ba \in \Z^n$, is 
squarefree if and only if $\ba$ is squarefree. 
Let $\grS$ (resp. $\SqS$) be the category of finitely generated 
$\Z^n$-graded $S$-modules (resp. squarefree $S$-modules). 
Let $\P$ be a $\Z^n$-graded minimal free resolution of $M \in \grS$. 
Then $M$ is squarefree if and only if each 
$P_i$ is a direct sum of copies of $S(-F)$ 
for various $F \subset [n]$. In the present case, 
an order ideal $\Delta$ of $\cL \, (\cong 2^{[n]})$ is 
essentially a simplicial complex, and the ring $K[\Delta]$ 
defined in the previous section is nothing other than 
the {\it Stanley-Reisner ring} (c.f. \cite{BH,St}) of $\Delta$.  

\medskip

Let $E = \bigwedge \< y_1, \ldots, y_n\>$ be the exterior algebra 
over $K$. Under the {\it Bernstein-Gel'fand-Gel'fand correspondence} 
(c.f. \cite{EFS}), $E$ is the counter part of $S$. 
We regard $E$ as a $\Z^n$-graded ring by $\deg y_i = \be_i = \deg x_i$ 
for each $i$. Clearly, any monomial ideal of $E$ is ``squarefree", 
and of the form $J_\Delta := ( \, \prod_{i \in F} y_i \mid F \subset [n], \, 
F \not \in \Delta \, )$ for a simplicial complex $\Delta \subset 2^{[n]}$. 
We say $K\< \Delta\> := E/J_\Delta$ is the 
{\it exterior face ring} of $\Delta$. 

Let $\grE$ be the category of finitely generated $\Z^n$-graded 
$E$-modules and degree preserving $E$-homomorphisms. 
Note that, for graded $E$-modules, we do not have to distinguish left modules 
from right ones. Hence $$\DE(-) := \bigoplus_{\ba \in \Z^n} 
\Hom_{\grE}(-, E(\ba))$$  
gives an exact contravariant functor from $\grE$ to itself 
satisfying $\DE \circ \DE = \Id$.

\begin{dfn}[R\"omer~\cite{R0}]\label{sqE}
We say $N \in \grE$ is 
{\it squarefree}, if $N = \bigoplus_{F \subset [n]} N_F$ 
(i.e., if $\ba \in \Z^n$ is not squarefree, then $N_\ba = 0$). 
\end{dfn}

An exterior face ring $K\< \Delta \>$ is a squarefree $E$-module. 
But, since a free module $E(\ba)$ is not squarefree 
for $\ba \ne 0$, the syzygies of a squarefree $E$-module are 
{\it not} squarefree. 
Let $\SqE$ be the full subcategory of $\grE$ consisting of 
squarefree modules. If $N$ is a squarefree $E$-module, 
then so is $\DE(N)$. That is, $\DE$ gives a contravariant functor 
from $\SqE$ to itself.  

We have functors 
$\cS: \SqE \to \SqS$ and $\cE: \SqS \to \SqE$ giving an 
equivalence $\SqS \cong \SqE$. Here 
$\cS(N)_F = N_F$ for $N \in \SqE$ and $F \subset [n]$, 
and the multiplication map $\cS(N)_F \ni z \mapsto x_i z \in 
\cS(N)_{F \cup \{ i \}}$ for $i \not \in F$ is given by 
$\cS(N)_F =N_F \ni z \mapsto (-1)^{\alpha(i, F)} y_i z 
\in N_{F \cup \{ i \}}= \cS(N)_{F \cup \{ i \}}$, where 
$\alpha(i,F) =  \# \{ \, j \in F \mid j < i \, \}$. 
For example. $\cS(K\< \Delta\>) \cong K[\Delta]$. 
See \cite{R0} for detail. 

Note that $\bA := \cS \circ \DE \circ \cE$ is an exact 
contravariant functor from $\SqS$ to itself satisfying 
$\bA \circ \bA = \Id$.  
It is easy to see that $\bA(K[F]) \cong S(-F^\cmp)$, where 
$F^\cmp := [n] \setminus F$. 
We also have  $\bA(K[\Delta]) \cong I_{\Delta^\vee}$, 
where $$\Delta^\vee := \{ \, F \subset [n] \mid F^\cmp \not \in \Delta 
\, \}$$ is the {\it Alexander dual} complex of $\Delta$. 
Since $\bA$ is exact, it exchanges a (minimal) free resolution with 
a (minimal) irreducible resolution.

\medskip 

Eisenbud et al.\ (\cite{Ei, EFS}) introduced the notion of the {\it linear strands} 
and the {\it linear part} of a minimal free resolution of a graded $S$-module. 
Let $\P : \cdots \to P_1 \to P_0 \to 0$ 
be a $\Z^n$-graded minimal $S$-free resolution of $M \in \grS$. 
We have natural numbers $\beta_{i,\ba}(M)$ 
for $i \in \NN$ and $\ba \in \Z^n$ such that 
$P_i = \bigoplus_{\ba \in \Z^n} S(-\ba)^{\beta_{i,\ba}(M)}$. 
We call $\beta_{i,\ba}(M)$ the {\it graded Betti numbers} of $M$.   
Set $|\ba| = \sum_{i=1}^na_i$ for $\ba = (a_1, \ldots, a_n) \in \NN^n$. 
For each $l \in \Z$, we define the $l$-{\it linear strand} 
$\lin_l(\P)$ of $\P$ as follows: 
The term $\lin_l(\P)_{i}$ of homological degree $i$  
is $$\bigoplus_{|\ba|=l+i} S(-\ba)^{\beta_{i, \ba}(M)},$$
which is a direct summand of $P_i$, 
and the differential $\lin_l(\P)_i \to \lin_l(\P)_{i-1}$ 
is the corresponding component of the differential 
$P_i \to P_{i-1}$ of $\P$.
By the minimality of $\P$, we can easily verify that
$\lin _l \ep{\P}$ are chain complexes (see also \cite[\S 7A]{Ei}). We call 
 $\lin(\P) := \bigoplus_{l \in \Z} \lin_l(\P)$ the 
{\it linear part} of $\P$. 
Note that the differential maps of $\lin(\P)$ are represented by 
matrices of linear forms. 
We call $$\ld_S(M):= \max \{ i \mid H_i(\lin(\P)) \ne 0 \}$$
the {\it linearity defect} of $M$. 

Sometimes, we regard $M \in \grS$ as a $\Z$-graded module by 
$M_j = \bigoplus_{|\ba|=j} M_\ba$. 
In this case, we set $\beta_{i,j}(M):= \bigoplus_{|\ba|=j} 
\beta_{i,\ba}(M)$. Then $\lin_l(\P)_i = S(-l-i)^{\beta_{i,l+i}(M)}$.

\begin{rem}
For $M \in \grS$, it is clear that $\ld_S(M) \leq \pd_S(M) \leq n$, 
and there are many examples attaining the equalities. 
In fact, $\ld_S(S/(x_1^2, \ldots, x_n^2)) = n$. 
But if $M \in \SqS$, then we always have $\ld_S(M) \leq n-1$. 
In fact, for a squarefree module $M$, 
$\pd_S(M) = n$, if and only if  $\depth_S M = 0$, if and only if 
$M \cong K \oplus M'$ for some $M' \in \SqS$. But $\ld_S(K) = 0$ and 
$\ld_S(M' \oplus K) = \ld_S(M')$. So we may assume that 
$\pd_S M' \leq n-1$.  
\end{rem}

\begin{prop}\label{ld_S}
Let $M \in \SqS$, and $\P$ its minimal graded free resolution. 
We have $$\max \{ \, i \mid H_i(\lin_l(\P))\not= 0\, \} = 
n-l - \depth_S ( \, \Ext_S^l( \bA(M), S) \, ), $$
and hence
$$\ld_S (M) = \max \{ \, i - \depth_S 
( \, \Ext_S^{n-i}( \bA(M), S) \, )  \mid  0 \leq i \leq n  \, \}.$$
\end{prop}

\begin{proof}
Note that  
$\I := \bA(\P)$ is a minimal irreducible resolution of $\bA(M)$. 
Moreover, we have $\bA(\lin_l(\P)) \cong \lin_{n-l}(\I)$. 
Since $\bA$ is exact, $$\max \{ \, i \mid H_i(\lin_l(\P)) \ne 0 
\, \} = \max  \{ \, i \mid H^i(\lin_{n-l}(\I)) \ne 0 \, \},$$ 
and hence 
\begin{equation}\label{ldirr = ldA}
\ld_S(M) = \ldirr_S(\bA(M)). 
\end{equation}
Hence the assertions follow from 
Corollary~\ref{ld for irr} (note that $S \cong \wS$ as underlying modules). 
\end{proof}

For $N \in \grE$, we have a $\Z^n$-graded minimal $E$-free resolution $\P$ of 
$N$. By the similar way to the  $S$-module case, we can define the linear 
part $\lin(\P)$ of $\P$, and set $\ld_E(N) := \max \{ \, i \mid 
H_i(\lin(\P)) \ne 0 \, \}$. (In \cite{R02,Y7}, $\ld_E(N)$ is 
denoted by $\operatorname{lpd}(N)$. 
``lpd" is an abbreviation for ``linear part dominate".)
In \cite[Theorem~3.1]{EFS}, 
Eisenbud et al. showed that $\ld_E (N) < \infty$ for all $N \in \grE$. 
Since $\pd_E (N) = \infty$ in most cases, this is a strong result. 
If $n \geq 2$, then we have 
$\sup \{ \, \ld_E(N) \mid  N \in \grE \, \} = \infty$.  In fact, 
since $E$ is selfinjective, we can take ``cosyzygies". 
But, if $N \in \SqE$, then $\ld_E(N)$ behaves quite nicely. 

\begin{thm}\label{S & E}
For $N \in \SqE$, we have $\ld_E(N) = \ld_S(\cS(N)) 
\leq n-1$. In particular, for a simplicial complex $\Delta \subset 2^{[n]}$, 
we have $\ld_E(K\< \Delta \>) = \ld_S(K[\Delta])$. 
\end{thm}

\begin{proof}
Using the Bernstein-Gel'fand-Gel'fand correspondence, 
the second author described $\ld_E (N)$ in  \cite[Lemma~4.12]{Y7}.  
This description is the first equality of the following computation, 
which proves the assertion. 
\begin{eqnarray*}
\ld_E (N) &=&  \max \{ \, i - \depth_S 
( \, \Ext_S^{n-i}( \cS \circ \DE(N), S) \, )  \mid  0 \leq i \leq n  \, \}  
\quad \text{(by \cite{Y7})}  \\ 
&=& \max \{ \, i - \depth_S 
( \, \Ext_S^{n-i}( \bA \circ \cS(N), S) \, )  \mid  0 \leq i \leq n  \, \}  
\quad \text{(see below)} \\
&=& \ld_S (\cS(N)) \quad \text{(by Proposition~\ref{ld_S}).}
\end{eqnarray*}
Here the second equality follows from the isomorphisms  
$\cS \circ \DE(N) \cong \cS \circ \DE \circ \cE \circ \cS (N) \cong 
\bA \circ \cS(N)$. 
\end{proof}

\begin{rem}
Herzog and R\"omer showed that $\ld_E(N) \leq \pd_S (\cS(N))$ 
for $N \in \SqE$ (\cite[Corollary~3.3.5]{R02}). 
Since $\ld_S(\cS(N)) \leq \pd_S (\cS(N))$ (the inequality 
is strict quite often), Theorem~\ref{S & E} refines 
their result. Our equality might follow from the argument in \cite{R02}, 
which constructs a minimal $E$-free resolution of $N$ from 
a minimal $S$-free resolution of $\cS(N)$. 
But it seems that certain amount of computation will be required. 
\end{rem}

Theorem \ref{S & E} suggests that we may set
$$
\ld (\Delta) := \ld_S (K[\Delta]) = \ld _E (K\< \Delta \>).
$$

\begin{thm}\label{ld topological}
If $I_\Delta \ne (0)$ (equivalently, $\Delta \ne 2^{[n]}$), 
then $\ld_S (I_\Delta)$ is a topological 
invariant of the geometric realization $|\Delta^\vee|$ of the Alexander dual 
$\Delta^\vee$ of $\Delta$. If $\Delta \ne 2^{T}$ for any $T \subset [n]$, then 
$\ld( \Delta )$ is also a topological invariant of $|\Delta^\vee|$ 
(also independent from the number $n= \dim S$).
\end{thm}

\begin{proof}
Since $\bA(I_\Delta) = K[\Delta^\vee]$ and $\Delta^\vee \ne \emptyset$, 
the first assertion follows from Theorem~\ref{topological} 
and the equality \eqref{ldirr = ldA} in the proof of 
Proposition~\ref{ld_S}. 

It is easy to see that $\Delta \ne 2^T$ for any $T$  
if and only if $\ld( \Delta ) \geq 1$.
If this is the case, $\ld ( \Delta) = \ld_S(I_\Delta)+1$, and  
the second assertion follows from the first.  
\end{proof}

\begin{rem}
(1) For the first statement of Theorem~\ref{ld topological}, 
the assumption that $I_\Delta \ne (0)$ is necessary. 
In fact, if $I_\Delta =(0)$, then $\Delta = 2^{[n]}$ and 
$\Delta^\vee = \emptyset$.  
On the other hand, if we set $\Gamma :=2^{[n]} \setminus [n]$, 
then $\Gamma^\vee = \{ \emptyset \}$ and $|\Gamma^\vee| = 
\emptyset = |\Delta^\vee|$. 
In view of Proposition~\ref{ld_S}, 
it might be natural to set $\ld_S(I_\Delta) = \ld_S( \, (0) \, )
= -\infty$.  But, $I_\Gamma = \wS$ and hence $\ld_S(I_\Gamma) 
= 0$. One might think it is better to set 
$\ld_S( \, (0) \, ) = 0$ to avoid the problem.
But this convention does not help so much,  
if we consider $K[\Delta]$ and $K[\Gamma]$.
In fact, $\ld_S(K[\Delta]) = \ld_S (S) = 0$ 
and $\ld_S(K[\Gamma]) = \ld_S(S/\wS) =1$. 

(2) Let us think about the second statement 
of the theorem.  
Even if we forget the assumption that $\Delta \ne 2^T$, 
$\ld(\Delta)$ is almost a topological invariant. 
Under the assumption that $I_\Delta \ne 0$, we have 
the following. 
\begin{itemize}
\item $\ld(\Delta) \leq 1$ if and only if 
$K[\Delta^\vee]$ is sequentially Cohen-Macaulay. 
Hence we can determine whether $\ld(\Delta) \leq 1$ 
from the topological space $|\Delta^\vee|$. 
\item $\ld(\Delta) = 0$, if and only if all facets of 
$\Delta^\vee$ have dimension $n-2$, if and only if 
$|\Delta^\vee|$ is Cohen-Macaulay and has dimension $n-2$.   
\end{itemize}
Hence, if we forget the number ``$n$", we can not determine 
whether $\ld(\Delta) = 0$ from $|\Delta^\vee|$.  
\end{rem}
%
%
\section{An upper bound of linearity defects.}
In the previous section, we have seen that $\ld_E(N) = \ld_S(\cS(N))$ for $N\in \SqE$, in particular $\ld _E \ep{\exfacering{\d}} = \ld _S \ep{\symfacering{\d}}$ for a simplicial complex $\d$. 
In this section, we will give an upper bound of them, and see that the bound is sharp.\\
\quad For $0 \not= N\in \grE$, regarding $N$ as a $\Z $-graded module,
we set $\indeg _E\ep{N} := \min\sbra{\, i \mid N_i \not= 0\, }$, which is called the {\it initial degree} of $N$,
and $\indeg _S\ep{M}$ is similarly defined as 
$\indeg _S\ep{M} := \min\sbra{\, i\mid M_i \not= 0\, }$ for $0\not= M\in \grS$.
If $\d \ne 2^{\verset }$ (equivalently $\stideali \d \not= 0$ or 
$\stidealj \d \not= 0$), then we have $\indeg _S\ep{\stideali{\d}} = \indeg _E\ep{\stidealj{\d}} 
= \min \sbra{\, \sharp F \mid F\subset \verset ,F\not\in \d \, }$,
where $\sharp F$ denotes the cardinal number of $F$. So we set
$$
\indeg \ep{\d} := \indeg _S\ep{\stideali{\d}} = \indeg _E\ep{\stidealj{\d}}.
$$
Since $\ld \ep{2^{\verset }} = \ld _S \ep{S} = \ld _E \ep{E} = 0$ holds, we henceforth exclude this trivial case;
we assume that $\d \not= 2^{\verset }$.\\
\indent We often make use of the following facts: 

\begin{lem}\label{sec:res lemma}
Let $0\not= M \in \grS$ and let $\P$ be a minimal graded free resolution of $M$. Then
\begin{enumerate}
\item $\lin _i \ep{\P} = 0$ for all $i < \indeg _S\ep{M}$, i.e., there are only $l$-linear strands with $l\ge \indeg _S\ep{M}$ in $\P${\rm ;}
\item $\lin _{\indeg _S\ep{M}} \ep{\P}$ is a subcomplex of $\P${\rm ;}
\item if $M\in \SqS$, then $\lin \ep{\P} = \dsum{}{0\le l \le n} \lin _l\ep{\P}$, and $\lin _l\ep{\P}_i = 0$ for all $i > n-l$ and all $0\le l \le n$,
where the subscript $i$ is a homological degree.
\end{enumerate}
\end{lem}
\begin{proof}
(1) and (2) are clear. (3) holds from the fact that $P_i \cong \dsum{}{F\subset \verset}S\ep{-F}^{\betti iF}$.
\end{proof}

\begin{thm}\label{sec:bound of ld}
For $0\not= N\in \SqE$, it follows that 
$$
\ld _E \ep{N}\le \max \sbra{ 0, n-\indeg _E\ep{N} -1}.
$$
\\
\quad By Theorem \ref{S & E}.
this is equivalent to say that for $M \in \SqS$,
$$
\ld _S \ep{M}\le \max \sbra{0, n-\indeg_S \ep{M} - 1}.
$$
\end{thm}  
\begin{proof}
It suffices to show the assertion for $M\in \SqS$. Set $\indeg _S\ep{M} = d$ and let $\P$ be a minimal graded free resolution of $M$.
The case $d=n$ is trivial by Lemma \ref{sec:res lemma} (1), (3). Assume that $d\le n-1$.
Observing that $\lin _l \ep{\P}_i = S \ep{-l-i}^{\betti i{i+l}}$, where $\betti i{i+l}$ are $\Z$-graded Betti numbers of $M$,
Lemma \ref{sec:res lemma} (1), (3) implies that
the last few steps of $\P$ are of the form
\begin{align*}
\begin{CD}
0@> >> S\ep{-n}^{\betti {n-d}n} @> >>S\ep{-n}^{\betti{n-d-1}n} \oplus S\ep{-n+1}^{\betti{n-d-1}{n-1}} @> >>\cdots .
\end{CD}
\end{align*}
Hence $\lin _d\ep{\P}_{n-d} =  S\ep{-n}^{\betti {n-d}n} = P_{n-d}$. Since $\lin _d\ep{\P}$ is a subcomplex of the acyclic complex $\P$ by Lemma \ref{sec:res lemma} (2),
we have $\H{n-d}{\lin _d\ep{\P}} = 0$, so that $\ld _S\ep{M} \le n-d -1$.
\end{proof}

Note that $\stidealj{\d}\in \SqE$ (resp. $\stideali{\d}\in \SqS$).
Since $\ld \ep{\d} \le \ld _E\ep{\stidealj{\d}} + 1$ 
(resp. $\ld \ep{\d} \le \ld _S\ep{\stideali{\d}} + 1$) holds, 
we have a bound for $\ld \ep{\d}$,
applying Theorem \ref{sec:bound of ld} to $\stidealj{\d}$ (resp. $\stideali{\d}$). 

\begin{cor}\label{sec:bound of d}
For a simplicial complex $\d$ on $\verset $, we have
$$
\ld \ep{\d} \le \max \sbra{1, n-\indeg \ep{\d}}.
$$
\end{cor}

Let $\d ,\Gamma$ be simplicial complexes on $\verset $. We denote $\d *\Gamma $ for the join
$$
\sbra{\, F\cup G \mid F\in \d , G\in \Gamma \, }
$$
of $\d$ and $\Gamma $, and for our convenience, set
$$
\ver \ep{\d} := \sbra{\ v\in \verset\mid \sbra{v} \in \d \ }.
$$

\begin{lem}\label{sec:lem of ld}
Let $\d$ be a simplicial complex on $\verset$. Assume that $\indeg \ep{\d} = 1$, or equivalently $\ver \ep{\d} \not= \verset$. Then we have
$$
\ld \ep{\d} = \ld \ep{\d *\sbra{v}}
$$
for $v \in \verset \setminus \ver \ep{\d}$.
\end{lem}

\begin{proof}
We may assume that $v=1$. Let $\P$ be a minimal graded free resolution of $\symfacering{\d *\sbra{1}}$ and $\koszul \ep{x_1}$ the Koszul complex
$$
\begin{CD}
0@> >>S\ep{-1} @> x_1 >> S @> >>0
\end{CD}
$$
with respect to $x_1$. Consider the mapping cone $\P \tensor{}{S} \koszul \ep{x_1}$ of the map $\P \ep{-1} \overset{x_1}{\longrightarrow} \P$. There is the short exact sequence
$$
\begin{CD}
0@> >>\P @> >>\P \tensor{}{S} \koszul \ep{x_1} @> >>\P\ep{-1}\bra{-1}@> >>0,
\end{CD}
$$
whence we have $\H i{\P\tensor{}{S} \koszul \ep{x_1}} = 0$ for all $i\ge2$ and the exact sequence
$$
\begin{CD}
0@> >>\H 1{\P\tensor{}{S} \koszul \ep{x_1}} @> >> \H 0{\P\ep{-1}} @> x_1>> \H 0{\P}.
\end{CD}
$$
But since $\H 0{\P} = \symfacering{\d*\sbra{1}}$ and $x_1$ is regular on it, we have
$\H 1{\P \tensor{}{S} \koszul \ep{x_1}} = 0$. Thus $\P \tensor{}{S} \koszul\ep{x_1}$ is acyclic and hence a minimal graded free resolution of $\symfacering{\d}$.
Note that $\lin \ep{\P\tensor{}{S} \koszul\ep{x_1}} = \lin \ep{\P} \tensor{}{S} \koszul\ep{x_1}$: in fact, we have
\begin{align*}
\lin _l\ep{\P \tensor{}{S} \koszul \ep{x_1}}_i &= \lin _l\ep{\P \tensor{}{S} S}_i \oplus \lin_l\ep{\P \bra{-1}\tensor{}{S} S\ep{-1}}_i \\
                                              &= \ep{\lin _l\ep{\P }_i \tensor{}{S} S} \oplus \ep{\lin _l\ep{\P}_{i-1} \tensor{}{S}S\ep{-1}} \\
                                              &= \ep{\lin _l\ep{\P } \tensor{}{S} \koszul \ep{x_1}}_i,
\end{align*}
where the subscripts $i$ denote homological degrees, and the differential map
$$
\lin _l\ep{\P\tensor{}{S} \koszul \ep{x_1}}_i \longrightarrow \lin _l\ep{\P\tensor{}{S} \koszul \ep{x_1}}_{i-1}
$$
is composed by $\partial ^{\< l\> }_i$, $-\partial ^{\< l\> }_{i-1}$, and the multiplication map by $x_1$,
where $\partial ^{\< l\> }_i$ (resp. $\partial ^{\< l\> }_{i-1}$) is the $i^{\text{th}}$(resp. $\ep{i-1}^{\text{st}}$) differential map of the $l$-linear strand of $\P$.
Hence there is the short exact sequence

\begin{align*}
\begin{CD}
0@> >>\lin \ep{\P} @> >> \lin \ep{\P \tensor{}{S} \koszul \ep{x_1}} @> >> \lin \ep{\P}\ep{-1}\bra{-1}@> >>0,
\end{CD}
\end{align*}
which yields that $\H i{\lin \ep{\P \tensor{}{S} \koszul \ep{x_1}}} = 0$ for all $i\ge \ld \ep{\d *\sbra{1}} +2$, and the exact sequence
\begin{multline*}
0 \longrightarrow \H {\ld \ep{\d *\sbra{1}} +1}{\lin \ep{\P \tensor{}{S} \koszul \ep{x_1}}} \longrightarrow \H {\ld \ep{\d *\sbra{1}}}{\lin \ep{\P}\ep{-1}} \\
\overset{x_1}{\longrightarrow} \H {\ld \ep{\d *\sbra{1}}}{\lin \ep{\P}} \longrightarrow \H {\ld \ep{\d *\sbra{1}}}{\lin \ep{\P \tensor{}{S} \koszul \ep{x_1}}}.
\end{multline*}
Since $x_1$ does not appear in any entry of the matrices representing the differentials of $\lin \ep{\P}$, it is regular on $\H {\bullet}{\lin \ep{\P}}$, and hence we have
$$
\H {\ld \ep{\d *\sbra{1}}+1}{\lin \ep{\P \tensor{}{S} \koszul \ep{x_1}}} = 0
$$
and 
$$
\H {\ld \ep{\d *\sbra{1}}}{\lin \ep{\P \tensor{}{S} \koszul \ep{x_1}}} \not= 0,
$$
since $\H {\ld \ep{\d *\sbra{1}}}{\lin \ep{\P}} \not= 0$. Therefore $\ld \ep{\d} = \ld \ep{\d *\sbra{1}}$.
\end{proof}

Let $\d$ be a simplicial complex on $\verset $. For $F\subset \verset $, we set
$$
\d _F := \sbra{\, G\in \d \mid G\subset F\, }.
$$
\quad The following fact, due to Hochster, is well known, but because of our frequent use, we mention it.
\begin{prop}[c.f. \cite{BH,St}]\label{sec:Hoc's formula a} 
For a simplicial complex $\d$ on $\verset $, we have
$$
\betti ij\ep{\symfacering{\d}}=\sum _{F\subset \verset, \sharp F =j}\dim _K\redhom{j-i-1}{\d _F}K ,
$$
where $\betti ij\ep{\symfacering{\d}}$ are the $\Z $-graded Betti numbers of $\symfacering{\d}$.
\end{prop}


Now we can give a new proof of \cite[Proposition 4.15]{Y7}, which is the latter part of the next result. 

\begin{prop}[cf. {\cite[Proposition 4.15]{Y7}}]\label{sec:indeg is 1}
Let $\d$ be a simplicial complex on $\verset $. If $\indeg \d = 1$, then we have
$$
\ld \ep{\d} \le \max \sbra{1, n-3}.
$$
Hence, for any $\d$, we have 
$$
\ld \ep{\d} \le \max \sbra{1, n-2}.
$$
\end{prop}

\begin{proof}
The second inequality follows from the first one and  Corollary \ref{sec:bound of d}. 
So it suffices to show the first.
We set $\mathcal{V} := \verset \setminus \ver \ep{\d}$. Our hypothesis $\indeg \d = 1$ implies that $\mathcal{V} \not= \varnothing$.
By Lemma \ref{sec:lem of ld}, the proof can be reduced to the case $\sharp \mathcal{V}=1$.
We may then assume that $\mathcal{V} = \sbra{1}$. Thus we have only to show that $\ld \ep{\d *\sbra{1}} \le \max \sbra{1,n-3}$. Since we have $\indeg \ep{\d *\sbra{1}} \ge 2$,
we may assume $n\ge 4$ by Corollary \ref{sec:bound of d}. The length of the $0$-linear strand of $\symfacering{\d *\sbra{1}}$ is $0$,
and hence we concentrate on the $l$-linear strands with $l\ge1$.
Let $\P$ be a minimal graded free resolution of $\symfacering{\d *\sbra{1}}$. Since, as is well known, the cone of a simplicial complex, i.e., the join with a point, is acyclic, we have
$$
\betti in \ep{\symfacering{\d *\sbra{1}}} = \dim _K\redhom{n-i-1}{\d *\sbra{1}}K = 0
$$
by Proposition \ref{sec:Hoc's formula a}. Thus $\lin _l\ep{\P}_{n-l} = 0$ for all $l\ge 1$.
Now applying the same argument as the last part of the proof of Theorem \ref{sec:bound of ld} (but we need to 
replace $n$ by $n-1$), we have
$$
\H{n-2}{\lin \ep{\P}} = 0,
$$
and so $\ld \ep{\d *\sbra{1}} \le n-3$.
\end{proof}



According to \cite[Proposition 4.14]{Y7}, we can construct a squarefree module $N\in \SqE$ with $\ld _E \ep{N} = \pd _S\ep{\cS\ep{N}} = n-1$.
By Theorems~\ref{S & E} and \ref{sec:bound of ld}, $M:=\cS\ep{N}$ satisfies 
that $\indeg _S\ep{M} = 0$ and $\ld _S\ep{M} = n-1$.
For $0\le i \le n-1$, let $\Omega _i\ep{M}$ be the $i^{\text{th}}$ 
syzygy of $M$. Then $\Omega _i\ep{M}$ is squarefree, and we have that $\ld _S\ep{\Omega _i\ep{M}} = \ld _S\ep{M} - i = n-i-1$
and $\indeg _S\ep{\Omega _i\ep{M}} \ge \indeg _S\ep{M} + i = i$.
Thus by Theorem \ref{sec:bound of ld}, we know that $\indeg _S\ep{\Omega _i\ep{M}} = i$ and $\ld _S\ep{\Omega _i\ep{M}} = n - \indeg _S\ep{\Omega _i\ep{M}} -1$. So the bound in Theorem \ref{sec:bound of ld}  is optimal.\\
\quad In the following, we will give an example of a simplicial complex $\d$ with $\ld \ep{\d} = n-\indeg \ep{\d}$ for $2\le \indeg \ep{\d} \le n-2$,
and so we know the bound in Proposition~\ref{sec:bound of d} is optimal if $\indeg \ep{\d} \ge 2$, that is, $\ver \ep{\d} = \verset$.  \\

Given a simplicial complex $\d$ on $\verset $, we denote $\sk {\d}i$ for the $i^{\text{th}}$ skeleton of $\d$, which is defined as
$$
\sk {\d}i := \sbra{\, F\in \d \mid \# F \le i+1\, }.
$$
\begin{exmp}\label{sec:sharp}
Set $\Sigma := 2^{\verset}$, and let $\Gamma $ be a simplicial complex on $\verset $ whose geometric realization $|\Gamma |$ is
homeomorphic to the $\ep{d-1}$-dimensional sphere with $2\le d<n-1$, which we denote by $S^{d-1}$. (For $m>d$ there exists a triangulation of $S^{d-1}$ with $m$ vertices.
See, for example, \cite[Proposition 5.2.10]{BH}).
Consider the simplicial complex $\d :=\Gamma \cup \sk{\Sigma }{d-2}$. We will verify that $\d$ is a desired complex, that is, $\ld \ep{\d} = n - \indeg \ep{\d}$.
For brief notation, we put $t:= \indeg \d$ and $l:=\ld \ep{\d }$. \\
\quad First, from our definition, it is clear that $t \ge d$. Thus it is enough to show that $n-d\le l$: in fact we have that
$l\le n-t \le n-d \le l$ by Corollary \ref{sec:bound of d}, and hence that $t=d$ and $l = n-d$. Our aim is to prove that
$$
\betti {n-d}n\ep{\symfacering{\d}} \not= 0 \quad \text{and} \quad \betti{n-d -1}{n-1}\ep{\symfacering{\d}} = 0,
$$
since, in this case, we have $\H {n-d}{\lin _d\ep{\P}}\not= 0$, and hence $n-d\le l$. \\
\quad Now, let $F\subset \verset$, and $\redcomp{\bullet }{\d _F}$, $\redcomp{\bullet }{\Gamma _F}$ be the augmented chain complexes of $\d _F$ and $\Gamma _F$, respectively. 
Since $\sk{\Sigma}{d-2}$ have no faces of dimension $\ge d-1$, we have $\redcomp{d-1}{\d _F} = \redcomp{d-1}{\Gamma _F}$ and hence $\redhom{d-1}{\d _F}K = \redhom{d-1}{\Gamma _F}K$.
On the other hand, our assumption that $|\Gamma |\approx S^{d-1}$ implies that $\Gamma $ is Gorenstein, and hence that
\begin{align*}
\redhom{d-1}{\Gamma _F}K = \begin{cases}
K & \text{if $F=\verset $;} \\
0 & \text{otherwise.}
\end{cases}
\end{align*} 
Therefore, by Proposition \ref{sec:Hoc's formula a}, we have that
\begin{align*}
\betti {n-d}n\ep{\symfacering{\d}} &= \dim _K\redhom{d-1}{\Gamma}K = 1\not= 0;\\
\betti{n-d -1}{n-1}\ep{\symfacering{\d}} &= \sum _{F\subset \verset, \sharp F=n-1}\dim _K\redhom{d-1}{\Gamma _F}K = 0.
\end{align*}
\end{exmp}
%
%
%
%
\section{A simplicial complex $\d$ with $\ld \ep{\d} = n-2$ is an $n$-gon}
Following the previous section, we assume that $\d \not= \verset$, throughout this section.
We say a simplicial complex on $\verset $ is an $n$-{\it gon} if its facets are $\sbra{1,2},\sbra{2,3},\cdots ,\sbra{n-1,n},$ and $\sbra{n,1}$ after a suitable permutation of vertices. 
Consider the simplicial complex $\d$ on $\verset$ given in Example \ref{sec:sharp}. If we set $d = 2$, then $\d$ is an $n$-gon. Thus if a simplicial complex $\d$ on $\verset $ is an $n$-gon, we have
$\ld \ep{\d} = n-2$. Actually, the inverse holds, that is, if $\ld \ep{\d} = n-2$ with $n\ge 4$, $\d$ is nothing but an $n$-gon.

\begin{thm} \label{sec:the case any dim}
Let $\d$ be a simplicial complex on $\verset$ with $n \ge 4$. Then $\ld \ep{\d}=n-2$ if and only if $\d$ is an $n$-gon.
\end{thm}

In the previous section, we introduced Hochster's formula (Proposition \ref{sec:Hoc's formula a}), but in this section, we need explicit correspondence
between $\bra{\tor{\bullet}S{\symfacering{\d}}K}_F$ and reduced cohomologies of $\d _F$, and so we will give it as follows. \\
\quad Set $V:=\angle{x_1,\dots ,x_n}=S_1$ and let $\lbull{\koszul}:=S\tensor{}{K} \bigwedge V$ be the Koszul complex of $S$ with respect to $x_1,\dots ,x_n$. Then we have
$$
\ebra{\tor iS{\symfacering{\d}}K}_F=\H i{\bra{\symfacering{\d}\tensor {}{S}\lbull{\koszul }}_F}=\H i{\ebra{\symfacering{\d}\tensor {}{K}\bigwedge V}_F}
$$
for $F\subset \verset$. Furthermore, the basis of the $K$-vector space $\ebra{\symfacering{\d}\tensor {}{K}\bigwedge V}_F$ is of the form $\bx ^G\tensor{}{}\wedge ^{F\setminus G}\bx$ with $G \in \d _F$,
where $\bx ^G = \prod _{i\in G}x_i$ and $\wedge ^{F\setminus G}\bx = x_{i_1} \wedge \cdots \wedge x_{i_k}$ for $\sbra{i_1,\cdots ,i_k} = F\setminus G$ with $i_1<\cdots <i_k$.
Thus the assignment
$$
\varphi ^i:\redcocomp{i-1}{\d _F}\ni e_G^{*} \longmapsto \ep{-1}^{\sign GF}\bx ^G\tensor{}{}\wedge ^{F\setminus G}\bx\in \ebra{\symfacering{\d}\tensor {}{K}\bigwedge V}_F
$$
with $G\in \d _F$ gives the isomorphism $\ubull{\varphi}:\redcocomp {\bullet}{\d _F}\ebra{-1} \longrightarrow \ebra{\symfacering{\d}\tensor {}{K}\bigwedge V}_F$ of chain complexes, where
$\redcocomp {i-1}{\d _F}$ (resp. $\redcomp{i-1}{\d _F}$) is the $\ep{i-1}^{\text{st}}$ term of the augmented cochain (resp. chain) complex of $\d _F$ over $K$, $e_G$ is
the basis element of $\redcomp{i-1}{\d _F}$ corresponding to $G$, and $e_G^{*}$ is the $K$-dual base of $e_G$. Here we set
$$
\alpha \ep{A,B}:=\sharp \sbra{\ep{a,b}\ |\ a>b, a\in A, b\in B}
$$
for $A,B \subset \verset$. Thus we have the isomorphism
\begin{equation}\label{eq:tor=redhom}
\bar{\varphi}:\redcohom{i-1}{\d _F}K \longrightarrow \ebra{\tor {\sharp F -i}S{\symfacering{\d}}K}_F.
\end{equation}

\begin{lem} \label{sec:position of ld}
Let $\d$ be a simplicial complex on $\verset$ with $\indeg \ep{\d} \ge 2$, and $\P$ a minimal graded free resolution of $\symfacering{\d}$. We denote $\lbull Q$ for the
subcomplex of $\P$ such that $Q_i := \bigoplus _{j \le i + 1} S\ep{-j}^{\betti ij} \subset \bigoplus _{j\in \Z}S\ep{-j}^{\betti ij} = P_i$.
Assume $n\ge 4$.
Then the following are equivalent.
\begin{enumerate}
\item $\ld \ep{\d} = n-2$;
\item $\H{n-2}{\lin _2\ep{\P}}\not= 0$;
\item $\H{n-3}{\lbull Q}\not= 0$.
\end{enumerate}
\end{lem}

In the case $n \ge 5$, the condition ($3$) is equivalent to $\H{n-3}{\lin _1\ep{\P}} \not= 0$.

\begin{proof}
Since $\indeg \ep{\d} \ge 2$, $\lin _0\ep{\P}_i = 0$ holds for $i \ge 1$. 
Clearly, $\H i{\lbull Q} = \H i{\lin _1\ep{\P}}$ for $i \geq 2$. 
Since $\lin _l\ep{\P}_i = 0$ for $i\ge n-2$ and $l\ge 3$ by Lemma \ref{sec:res lemma}
and that $\ld \ep{\d} \le n-2$ by Proposition~\ref{sec:indeg is 1}, 
it suffices to show the following.
\begin{align}\label{eq:what to show}
\H{n-2}{\lin_2\ep{\P}} \cong \H{n-3}{\lbull Q} \quad \text{and} \quad \H i{\lbull Q} = 0 \quad \text{for $i \ge n-2$}.
\end{align}
Since $\lbull Q$ is a subcomplex of $\P$, there exists the following short exact sequence of complexes.
$$
0 \longrightarrow \lbull Q \longrightarrow \P \longrightarrow \tilde{\P} := \P /\lbull Q \longrightarrow 0,
$$
which induces the exact sequence of homology groups
$$
\H i{\P} \longrightarrow \H i{\tilde{\P}} \longrightarrow \H{i-1}{\lbull Q} \longrightarrow \H{i-1}{\P}.
$$
Hence the acyclicity of $\P$ implies that $\H i{\tilde{\P}} \cong \H{i-1}{\lbull Q}$ for all $i \ge 2$.
Now $\H i{\tilde{\P}} = 0$ for $i \ge n-1$ by Lemma \ref{sec:res lemma} and the fact that $\tilde{P}_i = \oplus _{l \ge 2}\lin _l\ep{\P}_i$. So
the latter assertion of \eqref{eq:what to show} holds, since $n-2 \ge 2$. The former follows from the equality $\H{n-2}{\tilde{\P}} = \H{n-2}{\lin_2\ep{\P}}$,
which is a direct consequence of the fact that $\lin_2\ep{\P}$ is a subcomplex of $\tilde{\P}$, that $\tilde{P}_{n-2} = \lin _2\ep{\P}_{n-2}$, and that $\tilde{P}_{n-1} = 0$.
\end{proof}

Let $\d$ be a $1$-dimensional simplicial complex on $\verset$ (i.e., $\d$ is essentially a simple graph). A {\it cycle} $C$ in $\d$ of length $t$ ($\ge 3$) is
a sequence of edges of $\d$ of the form ($v_1$, $v_2$), ($v_2$, $v_3$)$,\dots ,$($v_t$, $v_1$) joining distinct vertices $v_1,\dots v_t$. \\

Now we are ready for the proof of Theorem \ref{sec:the case any dim}.

\renewcommand{\proofname}{\textit{Proof of Theorem \ref{sec:the case any dim}}}
\begin{proof}
The implication ``$\Leftarrow$" has been already done in the beginning of this section. So we shall show the inverse.
By Proposition \ref{sec:indeg is 1}, we may assume that $\indeg \ep{\d} \ge 2$.
Let $\P$ be a minimal graded free resolution of $\symfacering{\d}$ and $\lbull Q$ as in Lemma \ref{sec:position of ld}.
Note that $\lbull Q$ is determined only by $\ebra{\stideali{\d}}_2$ and that it follows $\ebra{\stideali{\d}}_2 = \ebra{\stideali{\sk{\d}1}}_2$.
If the $1$-skeleton $\sk{\d}1$ of $\d$ is an $n$-gon, then so is $\d$ itself. Thus by Lemma \ref{sec:position of ld}, we may assume that $\dim \d = 1$.
Since $\ld \ep{\d} = n-2$, by Lemma \ref{sec:position of ld} we have
$$
\redhom 1{\d}K \cong \redcohom 1{\d}K \cong \ebra{\tor{n-2}S{\symfacering{\d}}K}_{\verset} \not= 0,
$$
and hence $\d$ contains at least one cycle as a subcomplex. So it suffices to show that $\d$ has no cycles of length $\le n - 1$.
Suppose not, i.e., $\d$ has some cycles of length $\le n-1$. To give a contradiction, we shall show
\begin{align}\label{eq:contra}
0 \longrightarrow \lin _2\ep{\P}_{n-2} \longrightarrow \lin _{2}\ep{\P}_{n-3}
\end{align}
is exact; in fact it follows $\H{n-2}{\lin _2\ep{\P}} = 0$, which contradicts to Lemma \ref{sec:position of ld}.
For that, we need some observations
(this is a similar argument to that done in Theorem 4.1 of \cite{Y}). 
Consider the chain complex $\symfacering{\d}\tensor{}{K}\bigwedge V\tensor{}{K}S$
where $V$ is the $K$-vector space with the basis $x_1,\dots ,x_n$. 
We can define two differential map $\vartheta ,\partial $ on it as follows:
\begin{align*}
\vartheta \ep{f\tensor{}{}\wedge ^G\bx\tensor{}{} g} = \sum _{i\in G}\ep{-1}^{\sign iG}\ep{x_if\tensor{}{}\wedge ^{G\setminus \sbra{i}}\bx\tensor{}{}g}; \\
\partial  \ep{f\tensor{}{}\wedge ^G\bx\tensor{}{} g} = \sum _{i\in G}\ep{-1}^{\sign iG}\ep{f\tensor{}{}\wedge ^{G\setminus \sbra{i}}\bx\tensor{}{}x_ig}.
\end{align*}
By a routine, we have that $\partial \vartheta +\vartheta \partial = 0$, 
and easily we can check that the $i^{\text{th}}$ homology group of the chain complex
$\ep{\symfacering{\d}\tensor{}{K}\bigwedge V\tensor{}{K}S,\vartheta}$ 
is isomorphic to the $i^{\text{th}}$ graded free module of a minimal free resolution $\P$ 
of $\symfacering{\d}$.
Since, moreover, the differential maps of $\lin \ep{\P}$ is induced by $\partial$ 
due to Eisenbud-Goto \cite{EG} and Herzog-Simis-Vasconcelos \cite{HSV},
$\lin _l\ep{\P}_i \longrightarrow \lin _l\ep{\P}_{i-1}$ can be identified with
\begin{align*}
\dsum{}{F\subset \verset ,\sharp F=i+l}\bra{\tor iS{\symfacering{\d}}K}_F\tensor{}{K}S \overset{\bar{\partial}}{\longrightarrow } \dsum{}{F\subset \verset ,\sharp F =i-1+l}\bra{\tor {i-1}S{\symfacering{\d}}K}_F\tensor{}{K}S,
\end{align*}
where $\bar{\partial}$ is induced by $\partial $. 
In the sequel, $-\{i\}$ denotes the subset $\verset \setminus \{i\}$ of $\verset$.  
Then we may identify the sequence \eqref{eq:contra} with
$$
0 \longrightarrow \ebra{\tor {n-2}S{\symfacering{\d}}K}_{\verset}\tensor{}{K}S \overset{\bar{\partial}}{\longrightarrow }
\dsum{}{i \in \verset }\bra{\tor {n-3}S{\symfacering{\d}}K}_{-\{i\}}\tensor{}{K}S
$$
and hence, by the isomorphism \eref{tor=redhom}, with
\begin{align}\label{eq:lin_1}
\begin{CD}
0 \longrightarrow \redcohom 1{\d}K \tensor{}{K}S @> \bar{\varepsilon }>> {\displaystyle 
\dsum{}{i \in \verset}\redcohom 1{\d _{-\{i\}}}K \tensor{}{K}S}.
\end{CD}
\end{align}
Here $\bar{\varepsilon }$ is composed by  
$\bar{\varepsilon}_i: \redcohom 1{\d}K \tensor{}{K}S \to 
 \redcohom 1{\d_{-\{i\}}}K \tensor{}{K}S$ 
 which is induced by the chain map 
$$\varepsilon_i: 
\redcocomp \bullet{\d}\otimes _K S \longrightarrow 
\redcocomp \bullet{\d _{- \{i\}}}\otimes _K S,
$$
$$
\varepsilon _i\ep{e_G^{*}\otimes 1} =\begin{cases}
\ep{-1}^{\sign iG}e_G^{*} \otimes x_i & \text{if $i \not \in G$}; \\
0 & \text{otherwise}.
\end{cases}
$$
\quad Well, let $C$ be a cycle in $\d$ of the form ($v_1$, $v_2$), 
($v_2$, $v_3$)$,\dots ,$($v_t$, $v_1$) with distinct vertices $v_1, \cdots v_t$. 
We say $C$ has a {\it chord} if there exists an edge ($v_i$, $v_j$) of $G$
such that $j \not\equiv i+1$ ($\operatorname{mod}\ t$), and $C$ is said to be 
{\it minimal} if it has no chord. It is easy to see that the $1^{\text{st}}$ 
homology of $\d$ is generated by
those of minimal cycles contained in $\d$, that is, we have the surjective map:
$$
\dsum{}{\substack{C\subset \d \\ C : \text{minimal cycle}}} \redhom 1{C}K \longrightarrow \redhom 1{\d}K.
$$

 Now by our assumption that $\d$ contains a cycle of length $\le n-1$ 
(that is, $\d$ itself is not a minimal cycle), we have the surjective map
\begin{equation}\label{eta}
\begin{CD}
{\displaystyle \dsum{}{i \in  \verset}\redhom 1{\d _{- \{ i \}}}K} @> \bar{\eta} >>\redhom 1{\d}K
\end{CD}
\end{equation}
where $\bar{\eta}$ is induced by the chain map 
$\eta :\dsum{}{}\redcomp \bullet{\d _{- \{ i \}}} 
\longrightarrow \redcomp \bullet{\d}$,
and $\eta$ is the sum of 
$$\eta _i:\redcomp \bullet{\d _{- \{ i \}}} \ni e_G \mapsto 
\ep{-1}^{\sign iG}e_G \in \redcomp \bullet{\d}.$$  
Taking the $K$-dual of \eqref{eta}, we have the injective map
\begin{align*}
\begin{CD}
\redcohom 1{\d}K @>\bar{\eta}^{*}>> 
{\displaystyle \dsum{}{i \in \verset}\redcohom 1{\d _{- \{ i \}}}K},
\end{CD}
\end{align*}
where $\bar{\eta}^{*}$ is the $K$-dual map of $\bar{\eta}$, and  composed by the $K$-dual 
$$\bar{\eta}_i^*:\redcohom 1{\d}K \to \redcohom 1{\d_{-\{i\}}}K$$ 
of $\bar{\eta}_i$. 
Then for all $0\not= z \in \redcohom 1{\d}K$, we have $\bar{\eta}^{*}_i \ep{z} \not= 0$ for some $i$.
Recalling the map $\bar{\varepsilon }: \redcohom 1{\d}K \tensor{}{K}S \to 
\dsum{}{} \redcohom 1{\d _{-\{i\}}}K \tensor{}{K}S$ 
in \eqref{eq:lin_1} and its construction, we know for $z\in \redcohom 1{\d}K$,
$$
\bar{\varepsilon }\ep{z\otimes y} = \sum ^n_{i=1}\bar{\eta}_i^{*}\ep{z}\otimes x_iy,
$$
and hence $\bar{\varepsilon }$ is injective. 
\end{proof}

\begin{rem}
(1) If $\Delta$ is an $n$-gon, then $\Delta^\vee$ is an  
$(n-3)$-dimensional Buchsbaum complex with $\rH_{n-4}(\Delta^\vee;K) = K$.   
If $n=5$, then $\Delta^\vee$ is a triangulation of 
the M\"obius band. But, for $n \geq 6$, $\Delta^\vee$ is not 
a homology manifold.  In fact, let $\sbra{1,2},\sbra{2,3},\cdots ,\sbra{n-1,n}, \sbra{n,1}$ be the facets of $\d$, then if $F = [n] \setminus \{1,3,5 \}$,
easy computation shows that $\lk{\dual\Delta}F$ is a 0-dimensional  
complex with 3 vertices, and hence $\redhom 0{\lk{\dual\Delta}F}K
= K^2$.\\
\quad (2) If $\indeg \Delta \geq 3$, then the simplicial complexes  
given in Example \ref{sec:sharp} are not the only examples   
which attain the equality $\ld(\Delta) = n -\indeg \ep{\Delta}$.
We shall give two examples of such complexes. \\
\quad Let $\Delta$ be the triangulation of the real projective 
plane $\PP$ with 6 vertices which is given in 
\cite[figure 5.8, p.236]{BH}.
Since $\PP$ is a manifold, $K[\Delta]$ is Buchsbaum. 
Hence we have $H_\m^2(K[\Delta]) = [H_\m^2(K[\Delta])]_0 \cong 
\rH_1(\Delta;K)$. So, if $\chara(K) =2$, then we have  
$\depth_S (\Ext^4_S(K[\Delta], \omega_S))=0$. Note that we have $\d = \dual\d$ in this case.
Therefore, easy computation shows that 
$$
\ld \ep{\dual \d} = \ld(\Delta ) = 3 = 6-3 = 6 - \indeg \ep{\Delta}.
$$ 
\quad Next, as is well known, there is a triangulation of the torus with $7$ vertices. Let $\d$ be the triangulation.
Since $\dim \d = 2$, we have $\indeg \ep{\dual \d} = 7 - \dim \d - 1= 4$. Observing that
$\symfacering{\d}$ is Buchsbaum, we have, by easy computation, that
$$
\ld \ep{\dual\d} = 3 = 7 -4 = 7 - \indeg \ep{\dual\d}.
$$
Thus $\dual\d$ attains the equality, but is not a simplicial complex given in Example \ref{sec:sharp}, since it follows, from Alexander's duality, that
\begin{align*}
\dim _K \redhom i{\dual\d}K = \dim _K \redhom{4-i}{\d}K = \begin{cases}
2\not= 1 & \text{for $i=3$;} \\
0            & \text{for $i\ge 4$.}
\end{cases}
\end{align*}
\quad More generally, the dual complexes of $d$-dimensional Buchsbaum complexes $\d$ with $\redhom{d-1}{\d}K \not= 0$ satisfy the equality
$$
\ld\ep{\dual\d} = n - \indeg\ep{\dual\d},
$$
but many of them differ from the examples in Example \ref{sec:sharp}, and we can construct such complexes more easily as $\indeg \ep{\dual\d}$ is larger.
\end{rem}
%
%
%
%
%
%


\begin{thebibliography}{99}
\bibitem{AAH} A. Aramova, L. Avramov and J. Herzog, 
Resolutions of monomial ideals and cohomology over exterior algebras. 
Trans. Amer. Math. Soc. 352 (2000), 579-594.


\bibitem{BH} W. Bruns and J. Herzog,  Cohen-Macaulay rings, revised edition,  
Cambridge University Press, 1998.

\bibitem{Ei} D. Eisenbud, 
The geometry of syzygies: 
A second course in commutative algebra and algebraic geometry, 
Grad. Texts in Math., vol.229, Springer, 2005. 


\bibitem{EFS} D. Eisenbud, G. Fl\o ystad  and F.-O. Schreyer, 
Sheaf cohomology and free resolutions over exterior algebra, 
Trans. Amer. Math. Soc. 355 (2003), 4397-4426.  

\bibitem{EG} D. Eisenbud and S. Goto,
Linear free resolutions and minimal multiplicity,
J. Algebra 88 (1984), 89--133

\bibitem{HH} J. Herzog and T. Hibi,  Componentwise linear ideals, 
Nagoya Math. J. 153 (1999), 141--153.


\bibitem{HI} J. Herzog and S. Iyengar, Koszul modules, 
J. Pure Appl. Algebra 201(2005), 154--188.  

\bibitem{I} B. Iversen,  Cohomology of sheaves. Springer-Verlag, 1986. 

\bibitem{HSV} J. Herzog, A. Simis, and W. Vasconcelos,
Approximation complexes of blowing-up rings, II,
J. Algebra 82 (1983), 53-83.

\bibitem{MZ} R. Martinez-Villa and  D. Zacharia,  
Approximations with modules having linear resolutions, 
J. Algebra 266 (2003), 671--697.

\bibitem{Mil} E. Miller, Cohen-Macaulay quotients of 
normal affine semigroup rings via irreducible resolutions, Math. Res. Lett. 
9 (2002), 117-128.  

\bibitem{R0} T. R\"omer,  Generalized Alexander duality and applications, 
Osaka J. Math. 38 (2001), 469--485.  

\bibitem{R02} T. R\"omer, 
On minimal graded free resolutions, Thesis, University of Essen, 2001. 
 
\bibitem{Sc} P. Schenzel,
On the dimension filtration and Cohen-Macaulay filtered modules,
in: Commutative algebra and algebraic geometry,
ed. F. Van Oystaeyen,
Lecture Notes in Pure and Appl. Math., vol. 206, Dekker, 1999, pp. 245--264.
 
\bibitem{St} R. Stanley,  
Combinatorics and commutative algebra, 2nd ed. Birkh\"auser 1996.    

\bibitem{Y} K. Yanagawa, 
Alexander duality for Stanley-Reisner rings and squarefree 
$\NN^n$-graded modules, J. Algebra 225 (2000), 630--645.  

\bibitem{Y2} K. Yanagawa, 
Sheaves on finite posets and modules over normal semigroup rings, 
J. Pure and Appl. Algebra 161 (2001), 341--366. 

\bibitem{Y6} K. Yanagawa, 
Stanley-Reisner rings, sheaves, and Poincar\'e-Verdier duality,  
Math. Res. Lett. 10 (2003) 635--650. 

\bibitem{Y5} K. Yanagawa, 
Derived category of squarefree modules and local cohomology with monomial 
ideal support, J. Math. Soc. Japan 56 (2004) 289--308. 

\bibitem{Y7} K. Yanagawa, 
Castelnuovo-Mumford regularity for complexes and weakly Koszul modules, 
J. Pure and Appl. Algebra 207 (2006), 77--97.

\bibitem{Y8} K. Yanagawa, 
Notes on $C$-graded modules over an affine semigroup ring $K[C]$, Comm. Algebra, to appear.
\end{thebibliography}
\end{document}